\documentclass[12pt,english,twoside]{article}
\usepackage[english]{babel}
\usepackage{amsmath,amsfonts}
\usepackage{mathtools}
\usepackage[cp850]{inputenc}
\allowdisplaybreaks
\newcounter{theorem}[section]
\renewcommand{\thetheorem}{\thesection.\arabic{theorem}}
\newenvironment{lemma}[      1]{\refstepcounter{theorem} %
\bf \thetheorem\ Lemma#1.        \it}{}
\newenvironment{theorem}[    1]{\refstepcounter{theorem} %
\bf \thetheorem\ Theorem#1.      \it}{}
\newenvironment{proposition}[1]{\refstepcounter{theorem} %
\bf \thetheorem\ Proposition#1.  \it}{}

\newenvironment{example}[    1]{\refstepcounter{theorem} %
\bf \thetheorem\ Example#1.      \rm}{\par}
\newenvironment{examples}[   1]{\refstepcounter{theorem} %
\bf \thetheorem\ Examples#1.     \rm}{\par}
\newenvironment{remark}[  1]{\refstepcounter{theorem} %
\bf \thetheorem\ Remark#1.    \rm}{}

\newenvironment{proof}{ %
\it              Proof.          \rm}{\hfill $ \Box $}

\newcounter{abc}[theorem]

\newcounter{one}[theorem]
\newenvironment{onelist}{\begin{list}{%
\rm (\arabic{one}) \hfill           }{\usecounter{one} %
\topsep0mm \partopsep0mm \parsep0mm \itemsep0mm %
\leftmargin2em \labelwidth2em \labelsep0em}}{\end{list}}

\newenvironment{hylist}{\begin{list}{%
--               \hfill            }{%
\topsep0mm \partopsep0mm \parsep0mm \itemsep0mm %
\leftmargin2em \labelwidth2em \labelsep0em}}{\end{list}}

\newcounter{rom}
\newenvironment{romlist}{\begin{list}{%
\rm (\roman{rom})  \hfill           }{\usecounter{rom} %
\topsep0mm \partopsep0mm \parsep0mm \itemsep0mm %
\leftmargin2em \labelwidth2em \labelsep0em}}{\end{list}}

\newcommand{\fa}[3]{\{#1_#2\}_{#2\in#3}}

\newcommand{\N}{{\mathbb N}}
\newcommand{\R}{{\mathbb R}}

\newcommand{\CC}{{\cal C}}
\newcommand{\FF}{{\cal F}}

\newcommand{\MM}{{\cal M}}

\newcommand{\sss}{{\bf s}}

\newcommand{\uuu}{{\bf u}}

\newcommand{\xxx}{{\bf x}}

\newcommand{\UUU}{{\bf U}}

\newcommand{\XXX}{{\bf X}}

\newcommand{\one}{{\prime}}

\newcommand{\dps}{\displaystyle}

\newcommand{\leb}{{\boldsymbol\lambda}}

\newcommand{\zero}{{\bf 0}}
\newcommand{\eins}{{\bf 1}}

\newcommand{\leftmatrix}[1]{\left(\begin{array}{#1}}
\newcommand{\rightmatrix}{\end{array}\right)}

\begin{document}

\title{\Large\bf Estimators for a Class of Bivariate \\ Measures of Concordance for Copulas}
\author{Sebastian Fuchs and Klaus D. Schmidt}
\date{\normalsize Lehrstuhl f{\"u}r Versicherungsmathematik \\ Technische Universit{\"a}t Dresden}
\maketitle

\begin{abstract}
\noindent
In the present paper we propose and study estimators for a wide class of bivariate measures of concordance for copulas. 
These measures of concordance are generated by a copula and generalize Spearman's rho and Gini's gamma. 
In the case of Spearman's rho and Gini's gamma the estimators turn out to be the usual sample versions of these measures of concordance. 
\end{abstract}


\section{Introduction}
\label{intro}

The history of measures of concordance (or measures of association) 
starts with measures of concordance for a sample of bivariate random vectors. 
Later, 
related measures of concordance were introduced for bivariate distribution functions and copulas, 
and the sample versions for random vectors were interpreted as estimators of the population versions for distribution functions or copulas. 
Moreover, 
axioms for bivariate measures of concordance for copulas were developed, 
and most of these concepts have been extended to the multivariate case, 
with particular emphasis on Kendall's tau, Spearman's rho and Gini's gamma. 

\bigskip
In the present paper we propose and study estimators for a wide class of bivariate measures of concordance for copulas. 
These measures of concordance are generated by a copula and generalize Spearman's rho and Gini's gamma. 
In the case of Spearman's rho and Gini's gamma the estimators turn out to be the usual sample versions of these measures of concordance. 

\bigskip
This paper is organized as follows: 
In Section 
\ref{preliminaries} we resume some results on 
a group of transformations of copulas, 
invariance of copulas under a subgroup, 
measures of concordance for copulas which are defined in terms of the group, 
and 
a biconvex form for copulas. 
In Section 
\ref{moc} we consider a class of bivariate measures of concordance 
which are defined in terms of the biconvex form and are generated by a copula which is invariant under the full group of transformations; 
this class contains Spearman's rho and Gini's gamma as well as certain interpolations as special cases. 
In Section 
\ref{esti} we use the empirical copula to construct an estimator of the value of such a measure of concordance 
when the copula to be measured in unknown. 
To complete the discussion, 
we conclude with an Appendix on estimation under partial information on the copula to be measured: 
If the copula is known to be invariant under a specific subgroup of the group of transformations, 
then the value of every measure of concordance is equal to zero and the estimation problem is void. 

\bigskip
We denote 
by   $ \zero $   the vector in   $ \R^d $   whose coordinates are all equal to   $ 0 $   and 
by   $ \eins $   the vector in   $ \R^d $   whose coordinates are all equal to   $ 1 $.  
For a set   $ B\subseteq\R^d $, 
the indicator function   $ \chi_B : \R^d\to\{0,1\} $   is defined by   
$ \chi_B(\xxx):=1 $   if   $ \xxx \in B $   and 
$ \chi_B(\xxx):=0 $   else.


\section{Preliminaries}
\label{preliminaries}

In this section, 
we recall some definitions and results for the general dimension   $ d\geq2 $   
and point out the particularities in the case   $ d=2 $   
which are important for the subject of this paper.

\subsection*{A group of transformations of copulas}

Let   $ \CC $   denote the collection of all copulas   $ [\zero,\eins]\to[0,1] $. 
A map   $ \varphi: \CC\to\CC $   is said to be a 
\emph{transformation} on   $ \CC $. 
Let   $ \Phi $   denote the collection of all transformations on   $ \CC $   and define the 
\emph{composition}   $ \circ : \Phi\times\Phi\to\Phi $   by letting   $ (\varphi\circ\psi)(C) := \varphi(\psi(C)) $.   
The composition is associative and the transformation   $ \iota\in\Phi $   given by   $ \iota(C) := C $   
satisfies   $ \iota\circ\varphi = \varphi = \varphi\circ\iota $   for every   $ \varphi\in\Phi $   
and is therefore called the 
\emph{identity} on   $ \CC $.  
Thus, 
$ (\Phi,\circ) $   is a semigroup with neutral element   $ \iota $. 

\bigskip
We now introduce two types of elementary transformations. 
To this end, 
let   $ \MM $   denote the collection of all functions   $ [\zero,\eins]\to\R $. 
For   $ i,j\in\{1,...,d\} $   such that   $ i \neq j $   we define the 
\emph{transposition}   $ \pi_{i,j} : \CC\to\MM $   by letting
\begin{eqnarray*}
        (\pi_{i,j}(C))(\uuu)                                                             
& := &  \begin{cases}                                                                    
         C(u_1,\dots,u_{i-1},u_j,\dots,u_i,u_{j+1},\dots,u_d)  \quad\text{if   $ i<j $}  \\
         C(u_1,\dots,u_{j-1},u_i,\dots,u_j,u_{i+1},\dots,u_d)  \quad\text{if   $ i>j $}  
        \end{cases}                                                                      
\end{eqnarray*}
and for   $ k\in\{1,...,d\} $   we define the 
\emph{partial reflection}   $ \nu_k: \CC\to\MM $   by letting 
\begin{eqnarray*}
	    (\nu_{k}(C))(\uuu)                               
& := &  C(u_1,\dots,u_{k-1},1        ,u_{k+1}\dots,u_d)  
      - C(u_1,\dots,u_{k-1},1\!-\!u_k,u_{k+1}\dots,u_d)  
\end{eqnarray*}
Every transposition and 
every partial reflection is an involution in   $ \Phi $,  
and there exists a smallest subgroup   $ \Gamma $   of   $ \Phi $   containing 
all transpositions and 
all partial reflections. 
The group   $ \Gamma $   is a representation of the hyperoctahedral group with   $ d!\,2^{d} $   elements. 

\bigskip
A transformation is called 
\begin{hylist}
\item   a \emph{permutation} if it can be expressed as a finite composition of transpositions, 
and it is called 
\item   a \emph{reflection} if it can be expressed as a finite composition of partial reflections. 
\end{hylist}
Every transformation in   $ \Gamma $   can be expressed as a composition of a permutation and a reflection. 
We denote by 
\begin{hylist}
\item   $ \Gamma^\pi $   the set of all permutations and by 
\item   $ \Gamma^\nu $   the set of all reflections. 
\end{hylist}
Then   
$ \Gamma^\pi $   and   
$ \Gamma^\nu $   are subgroups of   $ \Gamma $, 
and   
$ \Gamma^\nu $   is commutative while   
$ \Gamma^\pi $   is not. 

\bigskip
Among the reflections, 
the 
\emph{total reflection} 
\begin{eqnarray*}
        \tau                        
& := &  \nu_d\circ\cdots\circ\nu_1  
\end{eqnarray*}
is of particular importance. 
The total reflection is an involution which transforms every copula into its survival copula. 
Define 
\begin{hylist}
\item   $ \Gamma^\tau := \{\iota,\tau\} $   and 
\item   $ \Gamma^{\pi,\tau} := \{ \gamma\in\Gamma \mid \gamma = \pi\circ\varphi \mbox{ for some } \pi\in\Gamma^\pi \mbox{ and some } \varphi\in\Gamma^\tau \} $
\end{hylist} \pagebreak
Then   $ \Gamma^\tau$   is the center of   $ \Gamma $, 
and   $ \Gamma^{\pi,\tau} $   is a subgroup of   $ \Gamma $. 
The total reflection also generates an order relation on   $ \CC $   
which compares not only two copulas but also their survival copulas:  
For   $ C,D\in\CC $   we write   $ C \preceq_\tau D $   if   $ C \leq D $   and   $ \tau(C)\leq\tau(D) $. 
Then   $ \preceq_\tau $   is an order relation  on   $ \CC $   which is called the 
\emph{concordance order}. 

\bigskip
\begin{remark}{\ (Bivariate case)}
Assume that   $ d=2 $   and let 
\begin{eqnarray*}
        \pi        
& := &  \pi_{1,2}  
\end{eqnarray*}
Then we have   $ \nu_2=\pi\circ\nu_1\circ\pi $   and  
\begin{eqnarray*}  
        \Gamma^\pi                       
&  = &  \{\iota,\pi\}                    \\
        \Gamma^{\pi,\tau}                
&  = &  \{\iota,\pi,\tau,\pi\circ\tau\}  \\*
        \Gamma^\nu                       
&  = &  \{\iota,\nu_1,\nu_2,\tau\}       
\end{eqnarray*}
Moreover, 
$ \Gamma $   is the smallest subgroup of   $ \Phi $   containing   $ \pi $   and   $ \nu_1 $, 
and the concordance order   $ \preceq_\tau $   coincides with the pointwise order   $ \leq $   on   $ \CC $. 
\end{remark}

\bigskip
Proofs and further details on the group of transformations may be found in 
Fuchs and Schmidt [2014] ($d=2$) and in 
Fuchs [2014]. 
We note in passing that in 
Fuchs and Schmidt [2014] the symbol   $ \nu $   is used instead of   $ \tau $.

\subsection*{Invariance of copulas with respect to a subgroup}

For a subgroup   $ \Lambda\subseteq\Gamma $, 
a copula   $ C $   is said to be $\Lambda$--\emph{invariant} if it satisfies   $ \gamma(C)=C $   for every   $ \gamma\in\Lambda $. 
For example, 
the product copula   $ \Pi $   is $\Gamma$--invariant, 
and 
the upper Fr{\'e}chet--Hoeffding bound   $ M $   is $\Gamma^{\pi,\tau}$--invariant. 
Moreover, 
for every copula   $ C $   and any subgroup   $ \Lambda\subseteq\Gamma $,  
the mean 
\begin{eqnarray*}
        C_\Lambda                                              
& := &  \frac{1}{|\Lambda|} \sum_{\gamma\in\Lambda} \gamma(C)  
\end{eqnarray*}
is a copula since   $ \CC $   is convex, 
and   $ C_\Lambda $   is $\Lambda$--invariant since   $ \Lambda $   is a subgroup. 
The collection of all $\Lambda$--invariant copulas is convex. 

\bigskip
\begin{remark}{\ (Bivariate case)}
Assume that   $ d=2 $. 
Then the lower Fr{\'e}chet--Hoeffding bound   $ W=\nu_1(M) $   is a copula which is $\Gamma^{\pi,\tau}$--invariant. 
\end{remark}

\bigskip
Proofs and further details on invariance of copulas may be found in  
Fuchs and Schmidt [2014] ($d=2$) and in 
Fuchs [2016b].

\subsection*{Measures of concordance for copulas}

A map   $ \kappa: \CC\to\R $   is said to be a 
\emph{measure of concordance} if it satisfies the following axioms:
\begin{romlist}
\item   $ \kappa[M] = 1 $.
\item   The identity   $ \kappa[\gamma(C)] = \kappa[C] $   holds for every   $ \gamma\in\Gamma^{\pi,\tau} $   and for all   $ C\in\CC $. 
\item   The identity   $ \sum_{\nu\in\Gamma^\nu} \kappa[\nu(C)] = 0 $   holds for all   $ C\in\CC $.  
\end{romlist}
These axioms are part of those proposed by 
Taylor [2007]. 
They imply that for every $\Gamma^\nu$--invariant copula   $ C $, 
and hence for every $\Gamma$--invariant copula 
and in particular the product copula, 
the identity   $ \kappa[C]=0 $   holds for every measure of concordance   $ \kappa $. 

\bigskip
A measure of concordance   $ \kappa : \CC\to\R $   is said to be 
\begin{hylist}
\item   \emph{convex} if   $ \kappa[(1\!-\!q)C + qD] = (1\!-\!q)\kappa[C] + q\kappa[D] $   holds for all   $ C,D\in\CC $   and all   $ q\in(0,1) $. 
\item   \emph{order preserving} if   $ \kappa[C]\leq\kappa[D] $   holds whenever   $ C,D\in\CC $   satisfy   $ C \leq D $. 
\item   \emph{concordance order preserving} if   $ \kappa[C]\leq\kappa[D] $   holds whenever   $ C,D\in\CC $   satisfy   $ C \preceq_\tau D $. 
\item   \emph{continuous} if   $ \lim_{n\to\infty} \kappa[C_n] = \kappa[C] $   
holds for any sequence   $ \{C_n\}_{n\in\N} \subseteq \CC $   and any copula   $ C\in\CC $   such that   $ \lim_{n\to\infty} C_n = C $   pointwise. 
\end{hylist}
If   $ \kappa $  is a concordance order preserving measure of concordance, 
then   $ \kappa[C] \leq \kappa[M] = 1 $   holds for all   $ C\in\CC $. 

\bigskip
\begin{remark}{\ (Bivariate case)} 
Assume that   $ d=2 $. 
Then a map   $ \kappa: \CC\to\R $   is a measure of concordance if and only if it has the following properties: 
\begin{romlist}
\item   $ \kappa[M] = 1 $.
\item   The identity   $ \kappa[  \pi(C)] =   \kappa[C] $   holds for all   $ C\in\CC $. 
\item   The identity   $ \kappa[\nu_1(C)] = - \kappa[C] $   holds for all   $ C\in\CC $. 
\end{romlist}
Moreover, 
a measure of concordance   $ \kappa $   is concordance order preserving if and only if it is order preserving, 
and in this case   $ -1 = \kappa[W] \leq \kappa[C] \leq \kappa[M] = 1 $   holds for all   $ C\in\CC $. 
\end{remark}

\bigskip
Proofs and further details on measures of concordance (as defined above) may be found in  
Fuchs and Schmidt [2014] ($d=2$) and in 
Fuchs [2016b].

\subsection*{A biconvex form for copulas}

Consider the map   $ [.\,,.]: \CC\times\CC\to\R $   given by 
\begin{eqnarray*}
        [C,D]                                     
& := &  \int_{[\zero,\eins]} C(\uuu)\,dQ^D(\uuu)  
\end{eqnarray*}
where   $ Q^D $   denotes the probability measure associated with the copula   $ D $; 
see Fuchs [2016a]. 
The map   $ [.\,,.] $   is in either argument linear with respect to convex combinations and is therefore called a 
\emph{biconvex form}. 
Moreover, 
the map   $ [.\,,.] $   is in either argument monotonically increasing with respect to the concordance order   $ \preceq_\tau $   
and it satisfies   $ 0\leq [C,D]\leq[M,M]=1/2 $   for all   $ C,D\in\CC $. 
Furthermore, 
there exist copulas   $ C,D\in\CC $   such that  $ [C,D]=0 $,   
and if   $ C,D\in\CC $   are $\Gamma^\nu$--invariant, 
then   $ [C,D] = 1/2^d = [\Pi,\Pi] $. 
Later on, 
we will also use the identities   $ [M,\Pi] = 1/(d+1) $   and  $ [M,M_{\Gamma^\nu}] = 1/4 + (1/2)^{d+1} $. 

\bigskip
\begin{remark}{\ (Bivariate case)}
The biconvex form   $ [.\,,.] $   is symmetric if and only if   $ d=2 $. 
\end{remark}

\bigskip
Proofs and further details on the biconvex form may be found in 
Fuchs [2016a].


\section{Measures of concordance generated by a copula}
\label{moc}

In the remainder of this paper 
we confine ourselves to the bivariate case   $ d=2 $   and 
we consider a wide class of measures of concordance which are defined in terms of the biconvex form. 

\bigskip
Consider a fixed copula   $ A\in\CC $. 
Then we have   $ [M,A] > [\Pi,A] $   such that the map   $ \kappa_A : \CC\to\R $   given by 
\begin{eqnarray*}
        \kappa_A[C]                              
& := &  \frac{[C,A] - [\Pi,A]}{[M,A] - [\Pi,A]}  
\end{eqnarray*}
is well--defined; 
see Fuchs [2016b]. 
We have the following result  
(see 
Edwards et al.\ [2005], 
Behboodian et al.\ [2005], 
Fuchs and Schmidt [2014], 
Fuchs [2016b]): 

\bigskip
\begin{proposition}{}
\label{char}
\begin{onelist}
\item   The map   $ \kappa_A $   is a measure of concordance if and only if   $ A $   is $\Gamma$--invariant. 
\item   If   $ A $   is $\Gamma$--invariant, 
then   $ \kappa_A $   is convex, order preserving and continuous. 
\item   If   $ A $   is $\Gamma$--invariant, 
then   $ [\Pi,A] = 1/4 $   and the identity 
$$  \kappa_A[C]  =  \frac{[C,A] - 1/4}{[M,A] - 1/4}  $$
holds for all   $ C\in\CC $. 
\end{onelist}
\end{proposition}

\bigskip
\begin{examples}{}
\label{example-1}
\begin{onelist}
\item   {\bf Spearman's rho:} 
The copula   $ \Pi $   is $\Gamma$--invariant and   $ \kappa_\Pi $   
satisfies 
\begin{eqnarray*}
        \kappa_\Pi[C]    
&  = &  12\,[C,\Pi] - 3  
\end{eqnarray*}
which means that   $ \kappa_\Pi $   is Spearman's rho; 
see 
Fuchs and Schmidt [2014] and 
Nelsen [2006; Subsection 5.1.2]. 
\item   {\bf Gini's gamma:} 
The copula   $ M_{\Gamma^\nu} = (1/2)(M\!+\!W)$   is $\Gamma$--invariant and   $ \kappa_{M_{\Gamma^\nu}} $   
satisfies 
\begin{eqnarray*}
        \kappa_{M_{\Gamma^\nu}}[C]  
&  = &  8\,[C,M_{\Gamma^\nu}] - 2   
\end{eqnarray*}
which means that   $ \kappa_{M_{\Gamma^\nu}} $   is Gini's gamma; 
see 
Fuchs and Schmidt [2014] and 
Nelsen [2006; Subsection 5.1.4]. 
\item   {\bf Linear interpolation:} 
For   $ q\in[0,1] $   define 
\begin{eqnarray*}
        E_q                             
& := &  (1\!-\!q)\Pi + qM_{\Gamma^\nu}  
\end{eqnarray*}
Then   $ E_q $   is a $\Gamma$--invariant copula and the measure of concordance   $ \kappa_{E_q} $   satisfies 
\begin{eqnarray*}
        \kappa_{E_q}[C]                                                                                                                        
&  = &  \frac{[C,E_q] - 1/4}{[M,E_q] - 1/4}                                                                                                    \\
&  = &  \frac{(1\!-\!q)([C,\Pi]-1/4)+q([C,M_{\Gamma^\nu}]-1/4)}{(1\!-\!q)([M,\Pi]-1/4)+q([M,M_{\Gamma^\nu}]-1/4)}                              \\
&  = &  \frac{(1\!-\!q)([M,\Pi]-1/4)                          }{(1\!-\!q)([M,\Pi]-1/4)+q([M,M_{\Gamma^\nu}]-1/4)}\,\kappa_\Pi[C]               \\*
&&\,+   \frac{                       q([M,M_{\Gamma^\nu}]-1/4)}{(1\!-\!q)([M,\Pi]-1/4)+q([M,M_{\Gamma^\nu}]-1/4)}\,\kappa_{M_{\Gamma^\nu}}[C]  \\*
&  = &  \frac{2(1\!-\!q)}{2+q}\,\kappa_\Pi[C] + \frac{3q}{2+q}\,\kappa_{M_{\Gamma^\nu}}[C]                                                     
\end{eqnarray*}
Therefore, 
$ \kappa_{E_q} $   is a weighted mean of Spearman's rho and Gini's gamma, 
and for   $ q\in(0,1) $   the respective weights are distinct from   $ 1-q $   and   $ q $. 
\end{onelist}
The last example can be extended to the case of arbitrary $\Gamma$--invariant copulas in the place of   $ \Pi $   and   $ M_{\Gamma^\nu} $. 
\end{examples}

\section{Estimation}
\label{esti}

We still assume that   $ d=2 $,   
and we also assume henceforth that the copula   $ A $   is $\Gamma$--invariant. 

\bigskip
For an arbitrary and unknown copula   $ C $, 
our aim is to construct an estimator of   $ \kappa_A[C] $   on the basis of a sample of bivariate random vectors with sample size   $ n\geq2 $. 
This estimator is based on an appropriate definition of the empirical copula, 
which in turn relies on the relative rank transform. 

\bigskip
The relative rank transform is constructed in three steps: 
\begin{hylist}
\item   Consider first the 
\emph{order transform}   $ T : \R^n\to\R^n $   which is defined coordinatewise by letting 
\begin{eqnarray*}
        (T(\xxx))_k                                                     
& := &  \min_{J\subseteq\{1,\dots,n\},\;|J|=k} \; \max_{j \in J}\; x_j  
\end{eqnarray*}
for all   $ k\in\{1,\dots,n\} $; 
see e.g.\ 
Fuchs and Schmidt [2016]. 
The map   $ T $   is measurable and for every   $ \xxx\in\R^n $   the coordinates of   $ T(\xxx) $   are increasing but need not be distinct. 
\item   Consider next the 
\emph{rank transform}   $ R : \R^n\to\{1,\dots,n\} $   which is defined coordinatewise 
as follows: 
Let   $ k_1 := \min\{ k\in\{1,\dots,n\} | x_k = (T(\xxx))_1 \} $   
and define 
\begin{eqnarray*}
        (R(\xxx))_{k_1}  
& := &  1                
\end{eqnarray*}
and for   $ i\in\!\{2,\dots,n\} $   
let   $ k_i := \min\{ k\in\!\{1,\dots,n\}\setminus\{k_1,\dots,k_{i-1}\} | x_k = (T(\xxx))_i \} $   
and define 
\begin{eqnarray*}
        (R(\xxx))_{k_i}  
& := &  i                
\end{eqnarray*}
The map   $ R $   is measurable and onto, 
and it satisfies   $ \sum_{k=1}^n (R(\xxx))_k = n(n\!+\!1)/2 $. 
\item   Consider finally the 
\emph{relative rank transform}   $ U : \R^n\to\{1/(n\!+\!1),\dots,n/(n\!+\!1)\} $   given by 
\begin{eqnarray*}
        U(\xxx)                 
& := &  \frac{1}{n+1}\,R(\xxx)  
\end{eqnarray*}
The map   $ U $   is measurable and onto, 
and it satisfies   $ \sum_{k=1}^n (U(\xxx))_k = n/2 $. 
\end{hylist}

\bigskip
Consider now 
a probability space   $ (\Omega,\FF,P) $   and 
an i.\,i.\,d.\ family   $ \fa{\XXX}{k}{\{1,\dots,n\}} $   of random vectors   $ \Omega\to\R^2 $   
such that   $ C $   is a copula for the distribution function of every   $ \XXX_k $. 
The family   $ \fa{\XXX}{k}{\{1,\dots,n\}} $   can be represented by the random matrix 
\begin{eqnarray*}
        \XXX_{(n)}                  
& := &  (\XXX_1,\dots,\XXX_n)       
\;\,=\;\,                           
        \leftmatrix{ccc}            
         X_{1,1} & \dots & X_{1,n}  \\
         X_{2,1} & \dots & X_{2,n}  
        \rightmatrix                
\;\,=\;\,                           
        \leftmatrix{c}              
         \XXX_{(n,1)}^\one          \\
         \XXX_{(n,2)}^\one          
        \rightmatrix                
\end{eqnarray*}
with   
$ (X_{1,k},X_{2,k}) := \XXX_k^\one $   and 
$ \XXX_{(n,i)} := ( X_{i,1},\dots,X_{i,n} )^\one $   for all   $ k\in\{1,\dots,n\} $   and   $ i\in\{1,2\} $. 
Define now 
\begin{eqnarray*}
        \UUU_{(n)}                  
& := &  \leftmatrix{c}              
         (U\circ\XXX_{(n,1)})^\one  \\
         (U\circ\XXX_{(n,2)})^\one  
        \rightmatrix                
\;\,=\;\,                           
        \leftmatrix{ccc}            
         U_{1,1} & \dots & U_{1,n}  \\
         U_{2,1} & \dots & U_{2,n}  
        \rightmatrix                
\;\,=\;\,                           
(\UUU_1,\dots,\UUU_n)               
\end{eqnarray*}
with   
$ (U_{i,1},\dots,U_{i,n}) := (U\circ\XXX_{(n,i)})^\one $   
and 
$ \UUU_k := (U_{1,k},U_{2,k})^\one $   for all   $ i\in\{1,2\} $   and   $ k\in\{1,\dots,n\} $. 
Then 
the rows of the random matrix   $ \UUU_{(n)} $   contain the relative ranks (without repetition) of 
the rows of the random matrix   $ \XXX_{(n)} $. 

\bigskip
The map   $ \widehat{C}_{(n)} : 
[\zero,\eins]
\times\Omega\to[0,1] $   given by 
\begin{eqnarray*}
        \widehat{C}_{(n)}(\uuu,\omega)                                  
& := &    \frac{1}{n} \sum_{k=1}^n \chi_{[\zero,\uuu]}(\UUU_k(\omega))  
\;\,=\;\, \frac{1}{n} \sum_{k=1}^n \chi_{[\UUU_k(\omega),\eins]}(\uuu)  
\end{eqnarray*}
is called the 
\emph{empirical copula} with sample size   $ n $   
(although it is not a copula since it fails to be continuous). 
This definition is appropriate for our purpose but it differs from that used in 
Nelsen [2006; Section 5.6]. 

\bigskip
\begin{lemma}{}
The empirical copula satisfies
$$  \int_{[\zero,\eins]} \widehat{C}_{(n)}(\uuu,\omega)\,dQ^{A}(\uuu)	 =  \frac{1}{n} \sum_{k=1}^n A(\UUU_k(\omega))  $$
for every   $ \omega\in\Omega $.
\end{lemma}

\bigskip
\begin{proof}{}
By continuity of   $ A $   and because of the identity   $ \sum_{k=1}^n U_{i,k}(\omega) = n/2 $   for   $ i\in\{1,2\} $, 
we obtain 
\begin{eqnarray*}
        \int_{[\zero,\eins]} \widehat{C}_{(n)}(\uuu,\omega)\,dQ^{A}(\uuu)                                                            
&  = &  \int_{[\zero,\eins]} \frac{1}{n} \sum_{k=1}^n \chi_{[\UUU_k(\omega),\eins]}(\uuu)\,dQ^{A}(\uuu)                              \\
&  = &  \frac{1}{n} \sum_{k=1}^n \int_{[\zero,\eins]} \chi_{(\UUU_k(\omega),\eins]}(\uuu)\,dQ^{A}(\uuu)                              \\
&  = &  \frac{1}{n} \sum_{k=1}^n Q^{A}[(\UUU_k(\omega),\eins]]  	                                                                 \\
&  = &  \frac{1}{n} \sum_{k=1}^n \Bigl( 1 - A(U_{1,k}(\omega),1) - A(1,U_{2,k}(\omega)) + A(U_{1,k}(\omega),U_{2,k}(\omega)) \Bigr)  \\
&  = &  \frac{1}{n} \sum_{k=1}^n \Bigl( 1 -   U_{1,k}(\omega)    -     U_{2,k}(\omega)  + A(U_{1,k}(\omega),U_{2,k}(\omega)) \Bigr)  \\*
&  = &  \frac{1}{n} \sum_{k=1}^n A(U_{1,k}(\omega),U_{2,k}(\omega))                                                                  
\end{eqnarray*}
as was to be shown. 
\end{proof}

\bigskip

Because of the previous result, 
we use the random variable 
\begin{eqnarray*}
        \langle C,A \rangle_{(n)}              
& := &  \frac{1}{n} \sum_{k=1}^n A\circ\UUU_k  
\end{eqnarray*}
as an estimator of   $ [C,A] $   when nothing is known about the copula   $ C $. 
By contrast, 
if it is known that   $ C=M $, 
then the coordinates of every   $ \XXX_k $   are comonotone and the relative ordinal ranks satisfy   $ U_{1,k}=U_{2,k}=k/(n\!+\!1) $   almost surely. 
Therefore, 
we use the real number 
\begin{eqnarray*}
        \langle M,A \rangle_{(n)}                                              
& := &  \frac{1}{n} \sum_{k=1}^n A\biggl( \frac{k}{n+1},\frac{k}{n+1} \biggr)  
\end{eqnarray*}
as an estimator of   $ [M,A] $.   
Correspondingly, 
we use the real number 
\begin{eqnarray*}
        \langle W,A \rangle_{(n)}                                                  
& := &  \frac{1}{n} \sum_{k=1}^n A\biggl( \frac{k}{n+1},\frac{n+1-k}{n+1} \biggr)  
\end{eqnarray*}
as an estimator of   $ [W,A] $. 

\bigskip
\begin{lemma}{}
\label{lemma}
\begin{onelist}
\item   $ \langle W,A \rangle_{(n)} \leq \langle C,A \rangle_{(n)} \leq \langle M,A \rangle_{(n)} $. 
\item   $ \langle M,A \rangle_{(n)} + \langle W,A \rangle_{(n)} = 1/2 $. 
\item   If   $ \langle M,A \rangle_{(n)} = 1/4, $ 
then   $ Q^A[(1/(n\!+\!1),n/(n\!+\!1)]^2] = 0 $. 
\item   If   $ Q^A[(1/(n\!+\!1),n/(n\!+\!1)]^2] > 0, $ 
then   $ \langle M,A \rangle_{(n)} > 1/4 $. 
\end{onelist}
\end{lemma}

\bigskip
\begin{proof}
To prove (1), 
consider a realization 
\begin{eqnarray*}
        (\uuu_1,\dots,\uuu_n)       
&  = &  \leftmatrix{ccc}            
         u_{1,1} & \dots & u_{1,n}  \\
         u_{2,1} & \dots & u_{2,n}  
        \rightmatrix                
\end{eqnarray*}
of the random matrix   $ (\UUU_1,\dots,\UUU_n) $. 
Then every row of the matrix 
$$  \leftmatrix{ccc}                
         u_{1,1} & \dots & u_{1,n}  \\
         u_{2,1} & \dots & u_{2,n}  
        \rightmatrix                
$$
contains each of the real numbers   $ 1/(n\!+\!1),\dots,n/(n\!+\!1) $   exactly once. 
Put 
\begin{eqnarray*}
        \uuu_k^{(0)}  
& := &  \uuu_k        
\end{eqnarray*}
and for   $ p\in\{1,\dots,n\} $   proceed as follows: 
Consider the unique   $ l_p\in\{1,\dots,n+1-p\} $   for which 
\begin{eqnarray*}
        \uuu_{l_p}^{(p-1)}                                    
&  = &  \biggl( \frac{i}{n+1},\frac{n+1-p}{n+1} \biggr)^\one  
\end{eqnarray*}   
holds for some   $ i\in\{1,\dots,n\} $. 
\begin{hylist}
\item   If   $ i=n+1-p $, 
put 
\begin{eqnarray*}
        \uuu_k^{(p)}    
& := &  \uuu_k^{(p-1)}  
\end{eqnarray*}
for every   $ k\in\{1,\dots,n\} $. 
\item   If   $ i \leq n-p $, 
consider the unique   $ m_p\in\{1,\dots,n\} $   for which 
\begin{eqnarray*}
        \uuu_{m_p}                                            
&  = &  \biggl( \frac{n+1-p}{n+1},\frac{j}{n+1} \biggr)^\one  
\end{eqnarray*}
holds for some   $ j\in\{1,\dots,n\} $    and put 
\begin{eqnarray*}
        \uuu_k^{(p)}                                                                 
& := &  \begin{cases}                                                                
\dps \biggl( \frac{n+1-p}{n+1},\frac{n+1-p}{n+1} \biggr)^\one & \text{if $ k=l_p $}  \\[2ex]
\dps \biggl( \frac{i    }{n+1},\frac{j    }{n+1} \biggr)^\one & \text{if $ k=m_p $}  \\[2ex]
\dps \uuu_k^{(p-1)}                                           & \text{else}          
        \end{cases}                                                                  
\end{eqnarray*}
\end{hylist}
In either case, 
we obtain 
\begin{eqnarray*}
        \uuu_{l_r}^{(p)}                                          
&  = &  \biggl( \frac{n+1-r}{n+1},\frac{n+1-r}{n+1} \biggr)^\one  
\end{eqnarray*}
for all   $ r\in\{1,\dots,p\} $, 
and since   $ A $   is $2$--increasing we also obtain 
\begin{eqnarray*}
        \frac{1}{n} \sum_{k=1}^n A(\uuu_k^{(p-1)})  
&\leq&  \frac{1}{n} \sum_{k=1}^n A(\uuu_k^{(p)})    
\end{eqnarray*}
After   $ n $   steps we thus obtain 
\begin{eqnarray*}
        \frac{1}{n} \sum_{k=1}^n A(\uuu_k)                                        
\;\,=\;\,  \frac{1}{n} \sum_{k=1}^n A(\uuu_k^{(0)})                               
&\leq&  \frac{1}{n} \sum_{k=1}^n A(\uuu_k^{(n)})                                  
\;\,=\;\,  \frac{1}{n} \sum_{k=1}^n A\biggl( \frac{k}{n+1},\frac{k}{n+1} \biggr)  
\end{eqnarray*}
and hence   $ \langle C,A \rangle_{(n)} \leq \langle M,A \rangle_{(n)} $. 
A similar algorithm yields   $ \langle W,A \rangle_{(n)} \leq \langle C,A \rangle_{(n)} $. 
This proves (1). 
\\
Since   $ A $   is $\Gamma$--invariant, 
we have   
\begin{eqnarray*}
        A(u,u)                    
&  = &     (\nu_2(A))(u,u)        
\;\,=\;\,  A(u,1) - A(u,1\!-\!u)  
\;\,=\;\,  u - A(u,1\!-\!u)       
\end{eqnarray*}
and hence 
\begin{eqnarray*}
        \langle M,A \rangle_{(n)}                                                                                  
&  = &  \frac{1}{n} \sum_{k=1}^n A\biggl( \frac{k}{n+1}\,,\frac{k}{n+1} \biggr)                                    \\
&  = &  \frac{1}{n} \sum_{k=1}^n \biggl( \frac{k}{n+1} - A\biggl( \frac{k}{n+1}\,,1-\frac{k}{n+1} \biggr) \biggr)  \\
&  = &  \frac{1}{2} - \frac{1}{n} \sum_{k=1}^n A\biggl( \frac{k}{n+1}\,,\frac{n+1-k}{n+1} \biggr)                  \\*
&  = &  \frac{1}{2} - \langle W,A \rangle_{(n)}                                                                    
\end{eqnarray*}
This proves (2). 
\\
Assume now that   $ \langle M,A \rangle_{(n)} = 1/4 $. 
Because of (2), 
this yields 
\begin{eqnarray*}
        0                                                                                                                     
&  = &  \langle M,A \rangle_{(n)} - \langle W,A \rangle_{(n)}                                                                 \\
&  = &  \frac{1}{n} \sum_{k=1}^n A\biggl( \frac{k}{n+1}\,,\frac{k}{n+1} \biggr)                                               
      - \frac{1}{n} \sum_{k=1}^n A\biggl( \frac{k}{n+1}\,,\frac{n+1-k}{n+1} \biggr)                                           \\
&  = &  \frac{1}{n} \sum_{k=1}^{\lfloor{n/2}\rfloor} \biggl(                                                                  
                 A\biggl( \frac{n+1-k}{n+1}\,,\frac{n+1-k}{n+1} \biggr) - A\biggl( \frac{n+1-k}{n+1}\,,\frac{k}{n+1} \biggr)  \\*
&&\!\qquad\quad- A\biggl( \frac{k}{n+1}\,,\frac{n+1-k}{n+1} \biggr) + A\biggl( \frac{k}{n+1}\,,\frac{k}{n+1} \biggr) \biggr)  
\end{eqnarray*}
Since   $ A $   is $2$--increasing, 
every term under the last sum is nonnegative and hence equal to   $ 0 $. 
This proves (3). 
\\
Assume finally that   $ Q^A[(1/(n\!+\!1),n/(n\!+\!1)]^2] > 0 $. 
Because of (3), 
this yields   $ \langle M,A \rangle_{(n)} \neq 1/4 $, 
and it then follows from (1) and (2) that   $ \langle M,A \rangle_{(n)} > 1/4 $. 
\end{proof}

\bigskip\pagebreak

Since the sequence   $ \{ (1/(n\!+\!1),n/(n\!+\!1)]^2 \}_{n\in\N} $   is increasing with union   $ (0,1)^2 $   and 
since   $ Q^A[(0,1)^2] = Q^A[[\zero,\eins]] = 1 $, 
we see that there exists some   $ n_A\in\N $   such that   $ Q^A[(1/(n\!+\!1),n/(n\!+\!1)]^2] > 0 $, 
and hence   $ \langle M,A \rangle_{(n)} > 1/4 $, 
holds for every   $ n\in\N $   with   $ n \geq n_A $. 
Since   $ Q^A[(1/2,1/2]^2] = Q^A[\emptyset] = 0 $, 
we have   $ n_A\geq2 $. 
The following example shows that   $ n_A $   may be greater than   $ 2 $: 

\bigskip
\begin{example}{}
According to Nelson [2006; Formula (3.1.5)], 
the map   $ E : [\zero,\eins]\to[0,1] $   given by 
\begin{eqnarray*}
        E(\uuu)                                                                          
& := &  \begin{cases}                                                                    
         \dps M(u_1,u_2)                     &  \text{if $ |u_1\!-\!u_2      |>1/2 $}  \\[.5ex]
         \dps W(u_1,u_2)                     &  \text{if $ |u_1\!+\!u_2\!-\!1|>1/2 $}  \\[.5ex]
         \dps \frac{u_1+u_2}{2}-\frac{1}{4}  &  \text{else}                              
        \end{cases}                                                                      
\end{eqnarray*}
is a copula, 
and it is straightforward to prove that   $ E $   is   $ \Gamma$--invariant and satisfies   $ \langle M,E \rangle_{(3)} = 1/4 $. 
Now Lemma 
\ref{lemma} yields   $ Q^E[(1/4,3/4]^2] = 0 $,   
and hence   $ n_E \geq 4 $. 
\end{example}

\bigskip
For the remainder of this section, 
we assume that the sample size   $ n $   satisfies   $ n \geq n_A $. 
Then we have   $ \langle M,A \rangle_{(n)} > 1/4 $   and the random variable 
\begin{eqnarray*}
        \widehat{\kappa_A[C]}_{(n)}                                              
& := &  \frac{\langle C,A \rangle_{(n)} - 1/4}{\langle M,A \rangle_{(n)} - 1/4}  
\end{eqnarray*}
is well--defined. 
We propose to use   $ \widehat{\kappa_A[C]}_{(n)} $   as an estimator of   $ \kappa_A[C] $. 

\bigskip
\begin{theorem}{}
The estimator   $ \widehat{\kappa_A[C]}_{(n)} $   satisfies 
$$  -1 = \widehat{\kappa_A[W]}_{(n)} \leq \widehat{\kappa_A[C]}_{(n)} \leq \widehat{\kappa_A[M]}_{(n)} = 1  $$%
\end{theorem}%

\bigskip
\begin{examples}{}
\begin{onelist}
\item   {\bf Spearman's rho:} 
Since 
\begin{eqnarray*}
        \langle C,\Pi \rangle_{(n)}                
&  = &  \frac{1}{n} \sum_{k=1}^n U_{1,k}\,U_{2,k}  
\end{eqnarray*}
and 
\begin{eqnarray*}
        \langle M,\Pi \rangle_{(n)}                               
&  = &  \frac{1}{n} \sum_{k=1}^n \biggl( \frac{k}{n+1} \biggr)^2  
\;\,=\;\,   \frac{2n+1}{6(n+1)}                                   
\end{eqnarray*}
the estimator of Spearman's rho satisfies 
\begin{eqnarray*}
        \widehat{\kappa_\Pi[C]}_{(n)}                        
&  = &  \frac{3}{n}\,\frac{n+1}{n-1}                         
        \Biggl( 4 \sum_{k=1}^n U_{1,k}\,U_{2,k} - n \Biggr)  
\end{eqnarray*}
Using the absolute ranks   $ R_{i,k} := (n\!+\!1)\,U_{i,k} $   instead of the relative ranks   $ U_{i,k} $,   
the previous identity can be written as 
\begin{eqnarray*}
        \widehat{\kappa_\Pi[C]}_{(n)}                                
&  = &  1 - \frac{6}{n(n^2\!-\!1)} \sum_{k=1}^n (R_{1,k}-R_{2,k})^2  
\end{eqnarray*}
which shows that the estimator is just the sample version of Spearman's rho; 
see also 
Kruskal [1958], 
Joe [1990] and 
P{\'e}rez and Prieto--Alaiz [2016]. 
\item   {\bf Gini's gamma:} 
Since 
\begin{eqnarray*}
        \langle C,M_{\Gamma^\nu} \rangle_{(n)}                                                                                                
&  = &  \biggl\langle C,\frac{1}{2}\,(M\!+\!W) \biggr\rangle_{\!\!(n)}                                                                        \\
&  = &  \frac{1}{2}\,\Bigl( \langle C,M \rangle_{(n)} + \langle C,W \rangle_{(n)} \Bigr)                                                      \\*
&  = &  \frac{1}{2}\,\Biggl( \frac{1}{n} \sum_{k=1}^n \min\{U_{1,k},U_{2,k}\} + \frac{1}{n} \sum_{k=1}^n \max\{U_{1,k}+U_{2,k}-1,0\} \Biggr)  
\end{eqnarray*}
and 
\begin{eqnarray*}
        \langle M,M_{\Gamma^\nu} \rangle_{(n)}                                                                                      
&  = &  \biggl\langle M,\frac{1}{2}\,(M\!+\!W) \biggr\rangle_{\!\!(n)}                                                              \\
&  = &  \frac{1}{2}\,\Bigl( \langle M,M \rangle_{(n)} + \langle M,W \rangle_{(n)} \Bigr)                                            \\
&  = &  \frac{1}{2}\,\Biggl( \frac{1}{2} + \frac{1}{n} \sum_{k=\lfloor(n+2)/2\rfloor}^n \biggl( \frac{2k}{n+1} - 1 \biggr) \Biggr)  \\*
&  = &  \frac{1}{4} + \frac{1}{4n(n\!+\!1)}\,\biggl\lfloor \frac{n^2}{2} \biggr\rfloor                                              
\end{eqnarray*}
the estimator of Gini's gamma satisfies 
\begin{eqnarray*}
  \widehat{\kappa_{M_{\Gamma^\nu}}[C]}_{(n)}                                                     
  &  = &  \frac{n+1}{\lfloor n^2/2 \rfloor}                     
          \Biggl( 2 \sum_{k=1}^n \min\{U_{1,k},U_{2,k}\} + 2 \sum_{k=1}^n \max\{U_{1,k}+U_{2,k}-1,0\} - n \Biggr)  
\end{eqnarray*}
Straightforward although slightly tedious calculation yields 
\begin{eqnarray*}
        \widehat{\kappa_{M_{\Gamma^\nu}}[C]}_{(n)}                                                                         
&  = &  \frac{1}{\lfloor n^2/2 \rfloor} \Biggl( \sum_{k=1}^n |R_{1,k}+R_{2,k}-1| - \sum_{k=1}^n |R_{1,k}-R_{2,k}| \Biggr)  
\end{eqnarray*}
which shows that the estimator is just the sample version of Gini's gamma; 
see 
Nelsen [2006; Subsection 5.1.4]. 
\item   {\bf Linear interpolation:} 
For   $ q\in[0,1] $   consider the $\Gamma$--invariant copula 
\begin{eqnarray*}
        E_q                               
&  = &  (1\!-\!q)\,\Pi + qM_{\Gamma^\nu}  
\end{eqnarray*}
introduced in Example 
\ref{example-1}(3). 
The estimator of   $ \kappa_{E_q}[C] $   satisfies 
\begin{eqnarray*}
        \widehat{\kappa_{E_q}[C]}_{(n)}                                                                                                                  
&  = &  \frac{\langle C,E_q \rangle_{(n)} - 1/4}{\langle M,E_q \rangle_{(n)} - 1/4}                                                                      \\
&  = &  \frac{(1\!-\!q)(\langle C,\Pi \rangle_{(n)}-1/4) + q(\langle C,M_{\Gamma^\nu} \rangle_{(n)} - 1/4)}                                              
             {(1\!-\!q)(\langle M,\Pi \rangle_{(n)}-1/4) + q(\langle M,M_{\Gamma^\nu} \rangle_{(n)} - 1/4)}                                              \\
&  = &  \frac{(1\!-\!q)(\langle M,\Pi \rangle_{(n)}-1/4)}                                                                                                
             {(1\!-\!q)(\langle M,\Pi \rangle_{(n)}-1/4) + q(\langle M,M_{\Gamma^\nu} \rangle_{(n)} - 1/4)}\,\widehat{\kappa_\Pi[C]}_{(n)}               \\*
&&\,   +\frac{                                             q(\langle M,M_{\Gamma^\nu} \rangle_{(n)} - 1/4)}                                              
             {(1\!-\!q)(\langle M,\Pi \rangle_{(n)}-1/4) + q(\langle M,M_{\Gamma^\nu} \rangle_{(n)} - 1/4)}\,\widehat{\kappa_{M_{\Gamma^\nu}}[C]}_{(n)}  
\end{eqnarray*}
and hence is a weighted mean of the estimators of Spearman's rho and Gini's gamma; 
for   $ q\in(0,1) $   the respective weights are distinct from   $ 1-q $   and   $ q $, 
due to the fact that   $ \langle M,\Pi \rangle_{(n)} \neq \langle M,M_{\Gamma^\nu} \rangle_{(n)} $. 
\end{onelist}
The last example can be extended to the case of arbitrary $\Gamma$--invariant copulas in the place of   $ \Pi $   and   $ M_{\Gamma^\nu} $. 
\end{examples}


\section{Appendix}

In certain cases in which some information on the copula   $ C $   is available, 
there is no need to estimate   $ \kappa[C] $   since this value is known. 
For example, 
the identities 
$ \kappa[W]=-1 $,   
$ \kappa[\Pi]=0 $   and 
$ \kappa[M]=1 $   hold for every measure of concordance   $ \kappa $. 
Moreover, 
if the copula   $ C $   is $\Gamma^\nu$--invariant, 
then the identity   $ \kappa[C]=0 $   holds for every measure of concordance   $ \kappa $   and the estimation problem for   $ \kappa[C] $   is void. 
The following result provides a class of continuous distribution functions for which the unique copula is $\Gamma^\nu$--invariant: 

\bigskip
\begin{theorem}{}
Assume that   $ C $   is a copula 
for which there exists 
a distribution function   $ F : \R^2\to[0,1] $   with marginal distribution functions   $ F_1,F_2 : \R\to[0,1] $   and 
a measurable function   $ f : \R^2\to\R_+ $   such that 
$$  f(x_1,x_2) = f(|x_1|,|x_2|)  $$
and 
$$  C(F_1(x_1),F_2(x_2)) = F(x_1,x_2) =  \int_{(-\infty,   x_1]\times(-\infty,x_2]} f(\sss)\,d\leb^2(\sss)  $$
$($with bivariate Lebesgue measure   $ \leb^2 )$   holds for every   $ \xxx\in\R^2 $. 
Then the copula   $ C $   is $\Gamma^\nu$--invariant and the identity   $ \kappa[C]=0 $   holds for every measure of concordance   $ \kappa $. 
\end{theorem}

\bigskip
\begin{proof}
Consider   $ \uuu\in[\zero,\eins] $   and any   $ \xxx\in\R^2 $   satisfying  $ F_i(x_i) = u_i $   for all   $ i\in\{1,2\} $. 
Then we have 
\begin{eqnarray*}
        (\nu_1(C))(\uuu)                                                   
&  = &  C(1,u_2) - C(1\!-\!u_1,u_2)                                        \\*
&  = &  C(1,F_2(x_2)) - C(1\!-\!F_1(x_1),F_2(x_2))                         \\
&  = &  C(1,F_2(x_2)) - C(F_1(-x_1),F_2(x_2))                              \\
&  = &  F_2(x_2) - F(-x_1,x_2)                                             \\
&  = &  \int_{(-\infty,\infty)\times(-\infty,x_2]} f(\sss)\,d\leb^2(\sss)  
      - \int_{(-\infty,  -x_1]\times(-\infty,x_2]} f(\sss)\,d\leb^2(\sss)  \\
&  = &  \int_{[-   x_1,\infty)\times(-\infty,x_2]} f(\sss)\,d\leb^2(\sss)  \\
&  = &  \int_{(-\infty,   x_1]\times(-\infty,x_2]} f(\sss)\,d\leb^2(\sss)  \\
&  = &  F(x_1,x_2)                                                         \\
&  = &  C(F_1(x_1),F_2(x_2))                                               \\*
&  = &  C(\uuu)                                                            
\end{eqnarray*}
This yields   $ \nu_1(C) = C $, 
and repeating the argument yields   $ \nu_2(C) = C $. 
Therefore, 
the copula   $ C $   is $\Gamma^\nu$--invariant. 
\end{proof}


%


\section*{References}

\small

Behboodian J, Dolati A, \'{U}beda--Flores M [2005]: 
Measures of association based on average quadrant dependence. 
Journal of Probability and Statistical Science 3, 171--173. 

\smallskip
Edwards HH, Mikusi\'{n}ski P, Taylor MD [2005]: 
Measures of concordance determined by $D_4$--invariant measures on $ (0,1)^2 $. 
Proceedings of the American Mathematical Society 133, 1505--1513. 

\smallskip
Fuchs S [2014]: 
Multivariate copulas: Transformations, symmetry, order and measures of concordance.
Kybernetika 50, 725--743. 

\smallskip
Fuchs S [2016a]: 
A biconvex form for copulas. 
Dependence Modeling 4, 63--75. 

\smallskip
Fuchs S [2016b]: 
Copula--induced measures of concordance. 
Dependence Modeling 4, 205--214. 

\smallskip
Fuchs S, Schmidt KD [2014]: 
Bivariate copulas: Transformations, asymmetry and measures of concordance. 
Kybernetika 50, 109--125. 

\smallskip
Fuchs S, Schmidt KD [2016]: 
On order statistics and their copulas. 
Statistics and Probability Letters 117, 165--172. 

\smallskip
Joe H [1990]: 
Multivariate concordance. 
Journal of Multivariate Analysis 35, 12--30. 

\smallskip
Kruskal WH [1958]: 
Ordinal Measures of Association. 
Journal of the American Statistical Association 53, 814--861. 

\smallskip
Nelson RB [1998]: 
Concordance and Gini's measure of association. 
Nonparametric Statistics 9, 227--238. 

\smallskip
Nelson RB [2006]: 
An Introduction to Copulas. Second Edition. 
Springer. 

\smallskip
P{\'e}rez A, Prieto--Alaiz M [2016]: 
A note on nonparametric estimation of copula--based multivariate extensions of Spearman's rho. 
Statistics and Probability Letters 112, 41--50. 

\smallskip
Scarsini M [1984]: 
On measures of concordance. 
Stochastica 8, 201--218. 

\smallskip
Taylor MD [2007]: 
Multivariate measures of concordance. 
Annals of the Institute of Statistical Mathematics 59, 789--806. 

\smallskip
Wolff EF [1980]: 
$n$--dimensional measures of dependence. 
Stochastica 4, 175--188.

\bigskip
\vfill\hspace*{\fill}\today
\end{document}